# MODERATE DEVIATIONS FOR PARTICLE FILTERING

By R. Douc, A. Guillin and J. Najim

*Ecole Polytechnique, Université Paris Dauphine and Ecole Nationale Supérieure des Télécommunications*

Consider the state space model $(X_t, Y_t)$, where $(X_t)$ is a Markov chain, and $(Y_t)$ are the observations. In order to solve the so-called filtering problem, one has to compute $\mathcal{L}(X_t|Y_1,\ldots,Y_t)$, the law of $X_t$ given the observations $(Y_1,\ldots,Y_t)$. The particle filtering method gives an approximation of the law $\mathcal{L}(X_t|Y_1,\ldots,Y_t)$ by an empirical measure $\frac{1}{n}\sum_1^n \delta_{x_{i,t}}$. In this paper we establish the moderate deviation principle for the empirical mean $\frac{1}{n}\sum_1^n \psi(x_{i,t})$ (centered and properly rescaled) when the number of particles grows to infinity, enhancing the central limit theorem. Several extensions and examples are also studied.

## 1. Introduction.

*The state space model.* Let $(X_t)$ be a $\mathbb{R}^d$-valued sequence of unobserved random variables and let $(Y_t)$ be the $\mathbb{R}^m$-valued observations, $t \geq 1$.

$$\begin{array}{ccccccccc}
X_0 & \to & X_1 & \to & \cdots & \to & X_{t-1} & \to & X_t & \to & X_{t+1} & \to & \cdots \\
& & \downarrow & & & & \downarrow & & \downarrow & & \downarrow & & \\
& & Y_1 & & & & Y_{t-1} & & Y_t & & Y_{t+1} & &
\end{array}$$

We endow $\mathbb{R}^d$ (resp. $\mathbb{R}^m$) with its Borel $\sigma$-field $\mathcal{B}(\mathbb{R}^d)$ [resp. $\mathcal{B}(\mathbb{R}^m)$] and we assume that $(X_t)_{t\in\mathbb{N}}$ is a Markov chain with initial distribution $\mathbb{P}(X_0 \in A) = \int_A a_0(x)\mu(dx)$ and transition kernels

$$\mathbb{P}(X_t \in A | X_{t-1} = x_{t-1}) = \int_A a_t(x_{t-1}, x_t)\,d\mu(x_t), \qquad t \geq 1,$$

where $\mu$ is a reference measure on $\mathbb{R}^d$ and $a_0$ (resp. $a_t$) is a probability density (resp. probability kernel density) with respect to $\mu$.

The observations $(Y_t)$ are $\mathbb{R}^m$-valued and independent conditionally on $(X_t)$. Formally, $Y_t$ depends on $X_t$ via the kernel

$$\mathbb{P}(Y_t \in B | X_t = x_t) = \int_B b_t(x_t, y_t) \, d\nu(y_t), \qquad t \geq 1,$$

where $\nu$ is a reference measure on $\mathbb{R}^m$ and $b_t$ is a probability density kernel with respect to $\nu$. It is worth noting that the general process $(X_t, Y_t)$ is a Markov chain.

These models, either called state space models or hidden Markov models (HMM), are widely used in engineering, biology, mathematical finance, geophysics, and so on. For an overview, see [18] and the references therein.

We will denote by $y_s^t$ the $\mathbb{R}^m$-valued series of observations $(y_s, \ldots, y_t)$ and when $s \leq t$, by $f_{t|s}(x_t|y_1^s)$ [or simply $f_{t|s}(x_t)$], the conditional density of $X_t$ given $Y_1^s = y_1^s$ with respect to $\mu$. In the case where $s = t$, $f_{t|t}$ is the filter density; in the case where $s = t-1$, $f_{t|t-1}$ is the one-step predictor. These quantities are related via the following relations: the *propagation* (or *prediction*) step,

$$(1.1) \qquad f_{t|t-1}(x_t|y_1^{t-1}) = \int f_{t-1|t-1}(x|y_1^{t-1}) a_t(x, x_t) \, d\mu(x),$$

and the *updating* step,

$$(1.2) \qquad f_{t|t}(x_t|y_1^t) = \frac{f_{t|t-1}(x_t|y_1^{t-1}) b_t(x_t, y_t)}{\int f_{t|t-1}(x|y_1^{t-1}) b_t(x, y_t) \mu(dx)}.$$

*Particle filtering.* The recursive computation of the filter density is a major issue. However, apart the very important Gaussian case for which filter density can be computed recursively with the Kalman–Bucy equations, there is no hope to get a closed-form formula for the filter density $f_{t|t}$ in the general case. Among the body of methods available to approximate the filter density (e.g., extended Kalman filter, approximate grid based filters, etc.), particle filtering (also known as recursive or sequential Monte Carlo filtering) has recently received a lot of attention. Let us mention the important contribution of Del Moral et al. [4, 5, 6, 7, 12] and the work of Künsch [18, 19]. The book edited by Doucet, de Freitas and Gordon [17] gives an overview of the subject and provides extra references.

In the sequel, we will say that $(x_1, \ldots, x_N)$ is a sample from $f \, d\mu$ if $(x_i; 1 \leq i \leq N)$ are independent and identically distributed (i.i.d.) with probability distribution $f \, d\mu$. We define recursively the approximate filter density:

At time $t = 0$, $(x_{i,0}; 1 \leq i \leq N)$ is a sample from $a_0 \, d\mu$ and

$$f_{1|1}^N(x) = \frac{b_1(x, y_1)(1/N) \sum_{i=1}^N a_1(x_{i,0}, x)}{\int_{\mathbb{R}^d} b_1(x, y_1)(1/N) \sum_{i=1}^N a_1(x_{i,0}, x) \mu(dx)}.$$

At time $t = T$, $(x_{i,T}; 1 \leq i \leq N)$ is a sample from $f_{T|T}^N d\mu$ and

$$f_{T+1|T+1}^N(x) = \frac{b_{T+1}(x, y_{T+1})(1/N)\sum_{i=1}^N a_{T+1}(x_{i,T}, x)}{\int_{\mathbb{R}^d} b_{T+1}(x, y_{T+1})(1/N)\sum_{i=1}^N a_{T+1}(x_{i,T}, x)\mu(dx)}.$$

As the number of particles $N$ grows to infinity, the empirical probability distribution $\frac{1}{N}\sum_1^N \delta_{x_{i,T}}$ converges to the filter probability distribution $f_{T|T} d\mu$. Among the main results for the particle filter, let us mention the law of large numbers [12], central limit theorems ([7, 9, 11], see also [19] for a nice exposition) and the large deviation principle [6].

*Links with genetic algorithms.* The approximate particle filter as expressed in the Introduction,

$$f_{T|T}^N(x) = \frac{b_T(x, y_T)(1/N)\sum_{i=1}^N a_T(x_{i,T-1}, x)}{\int_{\mathbb{R}^d} b_T(x, y_T)(1/N)\sum_{i=1}^N a_T(x_{i,T-1}, x)\mu(dx)},$$

can be interpreted as a genetic algorithm, a particle system approximation of the Feynman–Kac formula, as well as a so-called bootstrap filter in the filtering literature. This is of importance since up to some compatibility with the assumptions, we will then be able to rely on the important body of methods developed in the framework of particle systems approximation of the Feynman–Kac formulae (see [4, 6, 7, 8, 9, 10, 11, 12, 13]). Denote by

$$(1.3) \qquad \hat{b}_T(x_{T-1}) = \int_{\mathbb{R}^d} b_T(x, y_T) a_T(x_{T-1}, x)\mu(dx)$$

and

$$(1.4) \qquad \hat{a}_T(x_{T-1}, x_T) = \frac{b_T(x_T, y_T) a_T(x_{T-1}, x_T)}{\int_{\mathbb{R}^d} b_T(x_T, y_T) a_T(x_{T-1}, x_T)\mu(dx_T)}.$$

In this case, $f_{T|T}^N$ writes

$$f_{T|T}^N(x) = \sum_{i=1}^N \frac{\hat{b}_T(x_{i,T-1})}{\sum_{k=1}^N \hat{b}_T(x_{k,T-1})} \hat{a}_T(x_{i,T-1}, x)$$

and one can see the propagation and updating steps as a selection step followed by a mutation step:

$$(x_{i,T}, 1 \leq i \leq N) \overset{\text{selection}}{\longrightarrow} (\tilde{x}_{i,T}, 1 \leq i \leq N)$$

$$\overset{\text{mutation}}{\longrightarrow} (x_{i,T+1}, 1 \leq i \leq N).$$

The selection step consists in drawing a multinomial $\mathcal{M}(\omega_{1,T}, \ldots, \omega_{N,T})$, where $\omega_{i,T} = \frac{\hat{b}_T(x_{i,T-1})}{\sum_{i=1}^N \hat{b}_T(x_{i,T-1})}$ to choose accordingly the new particles $(\tilde{x}_{i,T}, 1 \leq i \leq N)$ among the generation $(x_{i,T}, 1 \leq i \leq N)$. At generation $T+1$, each particle $x_{i,T+1}$ is drawn independently according to the distribution $\hat{a}_T(\tilde{x}_{i,T}, \cdot) d\mu$.

*The main results.* In this paper we establish a moderate deviation principle (MDP) for the particle filter at time $T$ conditionally on the observations $(y_1, \ldots, y_T)$. Since the observations $(y_t)$ are frozen, this is a quenched result and we might sometimes drop the observations $(y_t)$ in our notation. In the sequel we will, therefore, denote by $\mathbb{P}_T$ the conditional probability $\mathbb{P}(\cdot | Y_1 = y_1, \ldots, Y_T = y_T)$.

The MDP complements the previously obtained CLTs [7, 9, 11, 19] and LDP [6] and is established for the rescaled and centered quantity

$$\mathbf{M}_N^T = \frac{1}{b_N \sqrt{N}} \sum_{i=1}^N \left( \psi_0(x_{i,0}) - \int \psi_0 a_0 \, d\mu, \ldots, \psi_T(x_{i,T}) - \int \psi_T f_{T|T} \, d\mu \right),$$

where the functions $(\psi_0, \ldots, \psi_T)$ are bounded and the speed $b_N^2$ is such that

$$\lim_{N \to \infty} b_N = \lim_{N \to \infty} \sqrt{N} b_N^{-1} = \infty.$$

The formal definition of an MDP states that there exists a good rate function $\mathbf{I}_T$ such that

$$-\inf_{\text{int}(\Gamma)} \mathbf{I}_T \leq \liminf_{N \to \infty} \frac{1}{b_N^2} \log \mathbb{P}_T \{ \mathbf{M}_N^T \in \Gamma \}$$

$$\leq \limsup_{N \to \infty} \frac{1}{b_N^2} \log \mathbb{P}_T \{ \mathbf{M}_N^T \in \Gamma \}$$

$$\leq -\inf_{\bar{\Gamma}} \mathbf{I}_T.$$

The set $\Gamma \subset \mathbb{R}^{T+1}$ is Borel, with interior $\text{int}(\Gamma)$ and closure $\bar{\Gamma}$. The rate function $\mathbf{I}_T$ depends on the asymptotic covariance matrix

$$\mathbf{V}_T(\psi_{0:T}) = (V_{s,t}(\psi_s, \psi_t))_{0 \leq s, t \leq T},$$

which appears in the central limit theorem (see details in Section 3.1). For applications of moderate deviations, see [1] or [14].

We then develop various extensions, such as the MDP for unbounded functions and a functional MDP for the particle density profile,

$$u \mapsto \frac{1}{b_N \sqrt{N}} \sum_{i=1}^{[Nu]} (\psi(x_{i,T}) - m_T(\psi)).$$

In this situation, the rate function turns out to be given by

$$J_T(f) = \int_0^1 \frac{\dot{f}^2(t)}{2\sigma_T^2(\psi)} \, dt + \frac{f^2(1)}{2} \left( \frac{1}{V_T(\psi)} - \frac{1}{\sigma_T^2(\psi)} \right),$$

where $f$ is absolutely continuous with $f(0) = 0$. The last part of the article is devoted to examples such as nonlinear observation models with additive noise and stochastic volatility models.

The paper is organized as follows. In Section 2 we give the notation and we state the main assumptions. In Section 3 we establish the MDP. Section 4 is devoted to various extensions of the MDP. Examples are studied in Section 5.

## 2. Notation, assumptions and a preliminary estimate.

2.1. *Notation and the main assumption.* We will endow $\mathbb{R}^d$ (resp. $\mathbb{R}^m$) with its Borel $\sigma$-field $\mathcal{B}(\mathbb{R}^d)$ [resp. $\mathcal{B}(\mathbb{R}^m)$]. Let $\mu$ be a reference measure on $(\mathbb{R}^d, \mathcal{B}(\mathbb{R}^d))$ and denote by $L^1(\mu)$ the set of $\mu$-integrable functions. Similarly, consider the reference measure $\nu$ defined on $(\mathbb{R}^m, \mathcal{B}(\mathbb{R}^m))$ and the function space $L^1(\nu)$. We will simply write $\mathcal{B}$ and $L^1$ whenever the context is clear. As usual, $\|\cdot\|_1$ denotes the $L^1$-norm and $\|\cdot\|_\infty$ the sup-norm ($\|f\|_\infty = \sup_x |f(x)|$). Recall that

$$\mathbb{P}(X_t \in A | X_{t-1} = x_{t-1}) = \int_A a_t(x_{t-1}, x_t) \, d\mu(x_t) \quad \text{and}$$

$$\mathbb{P}(Y_t \in B | X_t = x_t) = \int_B b_t(x_t, y_t) \, d\nu(y_t).$$

In the following, the sequence $(b_N)_{N \geq 1}$ will denote a sequence of nonnegative real numbers with the property that

$$\lim_{N \to \infty} b_N = \lim_{N \to \infty} \frac{\sqrt{N}}{b_N} = \infty.$$

We shall use the following notation (by convention $f_{0|0} = a_0$):

(2.1) $$m_t(\psi) = \int \psi f_{t|t} \, d\mu \quad \text{and} \quad m_{N,t}(\psi) = \int \psi f_{t|t}^N \, d\mu,$$

(2.2) $$L_t \psi(x) = \int a_t(x, u) b_t(u, y_t) \psi(u) \, \mu(du),$$

(2.3) $$\begin{aligned} M_N^t(\psi) &= \frac{1}{b_N \sqrt{N}} \sum_{i=1}^N (\psi(x_{i,t}) - m_t(\psi)) \\ &\triangleq Q_N^t(\psi) + R_N^t(\psi), \end{aligned}$$

where

(2.4) $$\begin{aligned} Q_N^t(\psi) &= \frac{1}{b_N \sqrt{N}} \sum_{i=1}^N (\psi(x_{i,t}) - m_{N,t}(\psi)) \quad \text{and} \\ R_N^t(\psi) &= \frac{\sqrt{N}}{b_N} (m_{N,t}(\psi) - m_t(\psi)). \end{aligned}$$

We might sometimes drop $\psi$ and simply write $Q_N^t$ and $R_N^t$. Denote by

$$\mathbf{M}_N^T(\psi_0,\ldots,\psi_T) = \mathbf{M}_N^T(\psi_{0:T}) \triangleq (M_N^0(\psi_0),\ldots,M_N^T(\psi_T)) \in \mathbb{R}^{T+1},$$

$$\boldsymbol{\lambda}_T = (\lambda_0,\ldots,\lambda_T) \in \mathbb{R}^{T+1}.$$

Recall that $\mathbb{P}_T = \mathbb{P}(\cdot|Y_1 = y_1,\ldots,Y_T = y_T)$. We will denote by $\mathbb{E}_T$ the expectation with respect to $\mathbb{P}_T$. Let us introduce now the main assumption on the model.

ASSUMPTION A-0. For every $t \geq 1$,

$$\gamma_t \triangleq \sup_x L_t \mathbb{1}(x) = \sup_x \int a_t(x,u) b_t(u,y_t) \mu(du) < \infty \quad \text{and}$$

$$L_t \mathbb{1}(x) > 0 \quad \forall x \in \mathbb{R}^d.$$

REMARK 2.1. Since $L_t \mathbb{1}(x) > 0$ under Assumption A-0, it is straightforward that

(2.5) $$\kappa_t \triangleq \int L_t \mathbb{1}(x) f_{t-1|t-1}(x) \mu(dx) > 0.$$

However, Assumption A-0 does not imply $\inf_x b_t(x,y_t) > 0$ as will be illustrated in the stochastic volatility model (see Section 5) where $\inf_x b_t(x,y_t) = 0$.

REMARK 2.2. A stronger version of Assumption A-0 is used in Section 4. See also Remark 4.1 for the link with genetic models.

Following [19], we define recursively the following variance-like quantities:

(2.6) $$\sigma_t^2(\psi) = \int (\psi - m_t(\psi))^2 f_{t|t}\, d\mu \quad \text{and}$$

$$\sigma_{N,t}^2(\psi) = \int (\psi - m_{N,t}(\psi))^2 f_{t|t}^N\, d\mu,$$

(2.7) $$V_t(\psi) = \sigma_t^2(\psi) + \frac{1}{\kappa_t^2} V_{t-1}(L_t \psi - m_t(\psi) L_t \mathbb{1}), \qquad V_0(\psi) = \sigma_0^2(\psi),$$

and the related covariance-like quantities:

(2.8) $$V_{r,t}(\psi_r, \psi_t) = \frac{1}{\kappa_t} V_{r,t-1}(\psi_r, L_t(\psi_t - m_t(\psi_t))) \qquad \text{for } r < t,$$

(2.9) $$V_{t,t}(\psi_t, \phi_t) = \tfrac{1}{2}(V_t(\psi_t + \phi_t) - V_t(\psi_t) - V_t(\phi_t)).$$

Of course, $V_{r,t}(\psi_r, \psi_t) = V_{t,r}(\psi_t, \psi_r)$. The covariance matrix is then defined by

$$\mathbf{V}_T(\psi_0,\ldots,\psi_T) = \mathbf{V}_T(\psi_{0:T}) = (V_{s,t}(\psi_s, \psi_t))_{0 \leq s,t \leq T}.$$

In the sequel, we will use $\langle \cdot, \cdot \rangle$ for the scalar product and "$\cdot$" for the matrix product.

2.2. *An exponential estimate.* In this section we prove an exponential estimate which will be useful to prove the MDP. This result is very close to Theorem 3.1 in [10] (see also Lemma 4 in [11] and Theorem 3.39 in [12]). However, since the model is slightly different, we provide a full proof.

LEMMA 2.1. *Assume that Assumption* A-0 *holds. Assume, moreover, that* $\psi:\mathbb{R}^d \to \mathbb{R}$ *is a bounded measurable function. Then, for every* $\varepsilon > 0$, *there exist* $\alpha(T) > 0$ *and* $\beta(T) > 0$ *such that*

$$(2.10) \quad \mathbb{P}_T\left(\left|\frac{1}{N}\sum_{i=1}^N \psi(x_{i,T}) - m_T(\psi)\right| > \varepsilon\right) \leq \alpha(T) e^{-N\varepsilon^2/(\beta(T)\|\psi\|_\infty^2)},$$

$$(2.11) \quad \mathbb{P}_T\left(\left|\int \psi(f_{T|T}^N - f_{T|T})\,d\mu\right| > \varepsilon\right) \leq \alpha(T) e^{-N\varepsilon^2/(\beta(T)\|\psi\|_\infty^2)}.$$

*Moreover, one can define recursively*

$$\begin{aligned}\alpha(0) &= 2, & & & \beta(0) &= 2, \\ \alpha(T) &= 4\max(1, \alpha(T-1)), & & \text{and} & \beta(T) &= \max\left(8, \frac{16\beta(T-1)\gamma_T^2}{\kappa_T^2}\right).\end{aligned}$$

PROOF. We shall prove (2.10) by induction. Recall that $\gamma_t$ is defined in Assumption A-0 and that $\kappa_t$ is defined in (2.5). At time $t = 0$, the result is a direct application of Hoeffding's inequality. Assume that (2.10) holds at time $t = T - 1$ and write

$$\frac{1}{N}\sum_{i=1}^N \psi(x_{i,T}) - m_T(\psi) = \frac{1}{N}\sum_{i=1}^N \psi(x_{i,T}) - m_{N,T}(\psi) + m_{N,T}(\psi) - m_T(\psi).$$

Thus,

$$\begin{aligned}\mathbb{P}_T&\left(\left|\frac{1}{N}\sum_{i=1}^N \psi(x_{i,T}) - m_T(\psi)\right| > \varepsilon\right) \\ &\leq \mathbb{P}_T\left(\left|\frac{1}{N}\sum_{i=1}^N \psi(x_{i,T}) - m_{N,T}(\psi)\right| > \frac{\varepsilon}{2}\right) \\ &\quad + \mathbb{P}_T\left(|m_{N,T}(\psi) - m_T(\psi)| > \frac{\varepsilon}{2}\right).\end{aligned}$$

Denote by $\mathcal{F}_N^T$ the $\sigma$-field generated by $(x_{i,t}; 0 \leq i \leq N, 0 \leq t \leq T)$. Conditionally on $\mathcal{F}_N^{T-1}$, the variables $(x_{i,T})$ are i.i.d. Therefore, Hoeffding's inequality yields

$$\mathbb{P}_T\left(\left|\frac{1}{N}\sum_{i=1}^N \psi(x_{i,T}) - m_{N,T}(\psi)\right| > \frac{\varepsilon}{2}\bigg|\mathcal{F}_N^T\right) \leq 2\exp\left(-\frac{N\varepsilon^2}{8\|\psi\|_\infty^2}\right),$$

which implies

$$(2.12) \quad \mathbb{P}_T\left(\left|\frac{1}{N}\sum_{i=1}^N \psi(x_{i,T}) - m_{N,T}(\psi)\right| > \frac{\varepsilon}{2}\right) \leq 2\exp\left(-\frac{N\varepsilon^2}{8\|\psi\|_\infty^2}\right).$$

Let us now deal with $m_{N,T}(\psi) - m_T(\psi)$. Apply the following identity: $\frac{A}{B} - \frac{A'}{B'} = \frac{A-A'}{B'} + \frac{A}{B}(\frac{B'-B}{B'})$ to $m_{N,T}(\psi) - m_T(\psi)$,

$$m_{N,T}(\psi) - m_T(\psi)$$
$$= \frac{(1/N)\sum_{i=1}^N L_T\psi(x_{i,T})}{(1/N)\sum_{i=1}^N L_T\mathbb{1}(x_{i,T})} - \frac{\int L_T\psi f_{T-1|T-1}\,d\mu}{\int L_T\mathbb{1} f_{T-1|T-1}\,d\mu}$$
$$= \frac{1}{\int L_T\mathbb{1} f_{T-1|T-1}\,d\mu}\left(\frac{1}{N}\sum_{i=1}^N L_T\psi(x_{i,T}) - \int L_T\psi f_{T-1|T-1}\,d\mu\right)$$
$$+ \frac{m_{N,T}(\psi)}{\int L_T\mathbb{1} f_{T-1|T-1}\,d\mu}\left(\frac{1}{N}\sum_{i=1}^N L_T\mathbb{1}(x_{i,T}) - \int L_T\mathbb{1} f_{T-1|T-1}\,d\mu\right).$$

Therefore,

$$\mathbb{P}_T\left(|m_{N,T}(\psi) - m_T(\psi)| > \frac{\varepsilon}{2}\right)$$
$$\leq \mathbb{P}_T\left\{\left|\frac{1}{N}\sum_{i=1}^N L_T\psi(x_{i,T}) - \int L_T\psi f_{T-1|T-1}\,d\mu\right| > \frac{\kappa_T\varepsilon}{4}\right\}$$
$$+ \mathbb{P}_T\left\{\left|\frac{1}{N}\sum_{i=1}^N L_T\mathbb{1}(x_{i,T}) - \int L_T\mathbb{1} f_{T-1|T-1}\,d\mu\right| > \frac{\kappa_T\varepsilon}{4\|\psi\|_\infty}\right\}.$$

As $\|L_T\psi\|_\infty^2 \leq \gamma_T^2\|\psi\|_\infty^2$, the induction assumption yields

$$\mathbb{P}_T\left(|m_{N,T}(\psi) - m_T(\psi)| > \frac{\varepsilon}{2}\right)$$
$$\leq \alpha(T-1)\exp\left(-\frac{N\kappa_T^2\varepsilon^2}{16\beta(T-1)\|L_T\psi\|_\infty^2}\right)$$
$$(2.13)$$
$$+ \alpha(T-1)\exp\left(-\frac{N\kappa_T^2\varepsilon^2}{16\beta(T-1)\|L_T\mathbb{1}\|_\infty^2\|\psi\|_\infty^2}\right)$$
$$\leq 2\alpha(T-1)\exp\left(-\frac{N\kappa_T^2\varepsilon^2}{16\beta(T-1)\gamma_T^2\|\psi\|_\infty^2}\right).$$

Inequality (2.10) is proved with the help of (2.12) and (2.13). Finally, (2.13) yields immediately (2.11). □

## 3. The moderate deviation principle.

3.1. *The MDP.* The moderate deviation principle is first proved for bounded test functions $\psi_0, \ldots, \psi_T$. The proof is simpler and one can focus on the main idea which is an induction approach. This technique has been used by Del Moral and Guionnet [6] for the LDP of the particle filter and by Künsch [19] for the CLT. The induction enables us to split $\mathbf{M}_N^T(\psi_0, \ldots, \psi_T)$ into one quantity depending on the last generation of particles $(x_{i,T})_{1 \leq i \leq N}$ and another one depending on all the other particles. These quantities turn out to be asymptotically independent. We relax the boundedness assumption over the test functions in Section 4.1.

THEOREM 3.1. *Assume that Assumption* A-0 *holds and let* $\psi_0, \ldots, \psi_T$ *be bounded measurable functions. The function defined by*

$$\mathbf{I}_T(\mathbf{x}_T) = \sup_{\boldsymbol{\lambda}_T \in \mathbb{R}^{T+1}} \left\{ \langle \mathbf{x}_T, \boldsymbol{\lambda}_T \rangle - \frac{\langle \boldsymbol{\lambda}_T, \mathbf{V}_T(\psi_{0:T}) \cdot \boldsymbol{\lambda}_T \rangle}{2} \right\}$$

*is a good rate function and the family of random variables* $(\mathbf{M}_N^T(\psi_{0:T}))_{N \geq 1}$ *satisfies the moderate deviation principle with speed* $b_N^2$ *and good rate function* $\mathbf{I}_T$, *that is,*

$$-\inf_{\mathrm{int}(\Gamma)} \mathbf{I}_T \leq \liminf_{N \to \infty} \frac{1}{b_N^2} \log \mathbb{P}_T\{\mathbf{M}_N^T(\psi_{0:T}) \in \Gamma\}$$

$$\leq \limsup_{N \to \infty} \frac{1}{b_N^2} \log \mathbb{P}_T\{\mathbf{M}_N^T(\psi_{0:T}) \in \Gamma\} \leq -\inf_{\bar{\Gamma}} \mathbf{I}_T,$$

*for* $\Gamma \in \mathcal{B}(\mathbb{R}^{T+1})$.

REMARK 3.1. If the covariance matrix $\mathbf{V}_T(\psi_{0:T})$ is invertible, then the rate function can be expressed as

$$\mathbf{I}_T(\mathbf{x}_T) = \frac{\langle \mathbf{x}_T, \mathbf{V}_T^{-1}(\psi_{0:T}) \cdot \mathbf{x}_T \rangle}{2}.$$

REMARK 3.2 (Particle profile). In the case where all the functions but $\psi_T$ are equal to zero, $\mathbf{M}_N^T(\psi_{0:T})$ reduces to the particle profile $M_N^T(\psi_T)$ and the rate function is given by the usual formula:

$$I_T(x) = \frac{x^2}{2V_T(\psi_T)}.$$

Moreover, one can prove under additional assumptions that the asymptotic variance $V_T(\psi_T)$ is uniformly bounded in time:

(3.1) $$V_t(\psi) \leq K\|\psi\|_\infty^2 \qquad \text{for all } t \geq 1.$$

This result is based on the property that the filter distribution forgets its initial condition (see [8], Theorem 3.1 and [13], Section 4.2.3 for the continuous time model). Equation (3.1) gives an MDP upper bound which does not depend on time.

REMARK 3.3 (Splitting the covariance matrix). Consider the covariance matrix $\mathbf{V}_T(\psi_{0:T})$. Denote by

$$\rho_{T-1}(x) = \left(\lambda_{T-1}\psi_{T-1} + \frac{\lambda_T}{\kappa_T}L_T(\psi_T - m_T(\psi_T))\right)(x)$$

and let $\tilde{\boldsymbol{\lambda}}_{T-1} = (\lambda_0, \ldots, \lambda_{T-2}, 1)$. Then the following identity holds true using (2.6)–(2.9):

$$\begin{aligned}(3.2) \quad & \langle \tilde{\boldsymbol{\lambda}}_{T-1}, \mathbf{V}_{T-1}(\psi_{0:T-2}, \rho_{T-1}) \cdot \tilde{\boldsymbol{\lambda}}_{T-1} \rangle \\ & = \langle \boldsymbol{\lambda}_T, \mathbf{V}_T(\psi_{0:T}) \cdot \boldsymbol{\lambda}_T \rangle - \lambda_T^2 \sigma_T^2(\psi_T).\end{aligned}$$

This identity will be useful in the sequel.

3.2. *Proof of Theorem* 3.1. We will proceed by induction. Let $u_0$ be a bounded function, then $\mathbf{M}_N^0(u_0)$ satisfies the MDP with good rate function $\mathbf{I}_0(x) = \frac{x^2}{2V_0(u_0)}$ since the particles $(x_{i,0})$ are i.i.d. with distribution $a_0 \, d\mu$. Assume that at time $T-1$, for every bounded functions $u_0, \ldots, u_{T-1}$, the random variables $(\mathbf{M}_N^{T-1}(u_{0:T-1}))_{N \geq 0}$ satisfy the MDP in $\mathbb{R}^T$ with good rate function $\mathbf{I}_{T-1}$. Consider now bounded functions $\psi_0, \ldots, \psi_T$ and the family of random variables $(\mathbf{M}_N^T(\psi_{0:T}))_{N \geq 1}$. The following lemma is crucial:

LEMMA 3.2. *Recall that by the definition of $f_{T|T}^N$ and by* (2.4),

$$R_N^T = \frac{\sqrt{N}}{b_N}(m_{N,T}(\psi_T) - m_T(\psi_T))$$

$$= \frac{\sqrt{N}}{b_N}\left(\frac{\sum_{i=1}^N L_T\psi_T(x_{i,T-1})}{\sum_{i=1}^N L_T\mathbb{1}(x_{i,T-1})} - m_T(\psi_T)\right)$$

*and let*

$$\tilde{R}_N^T = \frac{1}{b_N\sqrt{N}\kappa_T}\left(\sum_{i=1}^N L_T\psi_T(x_{i,T-1}) - m_T(\psi_T)\sum_{i=1}^N L_T\mathbb{1}(x_{i,T-1})\right),$$

*where $\kappa_T$ is defined in* (2.5). *Then the random variables $R_N^T$ and $\tilde{R}_N^T$ are exponentially equivalent up to the speed $b_N^2$. Otherwise stated,*

$$\limsup_{N \to \infty} \frac{1}{b_N^2} \log \mathbb{P}_T\{|R_N^T - \tilde{R}_N^T| > \delta\} = -\infty \qquad \text{for all } \delta > 0.$$

Proof of Lemma 3.2 is postponed to Appendix A. It is an alternative to the delta-method used for the CLT in [19].

PROOF OF THEOREM 3.1 (Continued). The random variables $M_N^T = Q_N^T + R_N^T$ and $\tilde{M}_N^T \triangleq Q_N^T + \tilde{R}_N^T$ are exponentially equivalent by Lemma 3.2. Therefore, it is sufficient by Gärtner–Ellis' theorem ([14], Theorem 2.3.6) to prove that

$$\lim_{N \to \infty} \frac{1}{b_N^2} \log \mathbb{E}_T \exp \left\{ b_N^2 \sum_{t=0}^{T-1} \lambda_t M_N^t(\psi_t) + b_N^2 \lambda_T \tilde{M}_N^T(\psi_T) \right\}$$
(3.3)
$$= \frac{\langle \boldsymbol{\lambda}_T, \mathbf{V}_T(\psi_{0:T}) \cdot \boldsymbol{\lambda}_T \rangle}{2}.$$

By (3.2),

$$\frac{1}{b_N^2} \log \mathbb{E}_T \left( \exp \left\{ b_N^2 \sum_{t=0}^{T-1} \lambda_t M_N^t(\psi_t) + b_N^2 \lambda_T \tilde{M}_N^T(\psi_T) \right\} \right)$$
$$- \frac{\langle \boldsymbol{\lambda}_T, \mathbf{V}_T(\psi_{0:T}) \cdot \boldsymbol{\lambda}_T \rangle}{2}$$
$$= \frac{1}{b_N^2} \log \mathbb{E}_T \left( \exp \left\{ \lambda_T b_N^2 Q_N^T - \frac{\lambda_T^2 b_N^2 \sigma_{N,T}^2(\psi_T)}{2} + \frac{\lambda_T^2 b_N^2}{2} \Delta_N \right. \right.$$
$$\left. \left. + b_N^2 \lambda_T \tilde{R}_N^T + b_N^2 \sum_{t=0}^{T-1} \lambda_t M_N^t(\psi_t) \right\} \right)$$
$$- \frac{\langle \tilde{\boldsymbol{\lambda}}_{T-1}, \mathbf{V}_{T-1}(\psi_{0:T-2}, \rho_{T-1}) \cdot \tilde{\boldsymbol{\lambda}}_{T-1} \rangle}{2},$$

where $\Delta_N = \sigma_{N,T}^2(\psi_T) - \sigma_T^2(\psi_T)$. Recall that $\mathcal{F}_N^T$ is the $\sigma$-field generated by $(x_{i,t}; 0 \leq i \leq N, 0 \leq t \leq T)$. In this case, $\Delta_N$ and $\tilde{R}_N^T$ are measurable with respect to $\mathcal{F}_N^{T-1}$. Thus, we get

$$\frac{1}{b_N^2} \log \mathbb{E}_T \left( \exp \left\{ b_N^2 \sum_{t=0}^{T-1} \lambda_t M_N^t(\psi_t) + b_N^2 \lambda_T \tilde{M}_T(\psi_T) \right\} \right)$$
$$- \frac{\langle \boldsymbol{\lambda}_T, \mathbf{V}_T(\psi_{0:T}) \cdot \boldsymbol{\lambda}_T \rangle}{2}$$
$$= \frac{1}{b_N^2} \log \mathbb{E}_T \left[ \mathbb{E}_T \left( \exp \left\{ \lambda_T b_N^2 Q_N^T - \frac{\lambda_T^2 b_N^2 \sigma_{N,T}^2(\psi_T)}{2} \right\} \Big| \mathcal{F}_N^{T-1} \right) \right.$$
$$\left. \times \exp \left\{ \frac{\lambda_T^2 b_N^2}{2} \Delta_N + b_N^2 \lambda_T \tilde{R}_N^T + b_N^2 \sum_{t=0}^{T-1} \lambda_t M_N^t(\psi_t) \right\} \right]$$

$$-\frac{\langle \tilde{\boldsymbol{\lambda}}_{T-1}, \mathbf{V}_{T-1}(\psi_{0:T-2}, \rho_{T-1}) \cdot \tilde{\boldsymbol{\lambda}}_{T-1}\rangle}{2}.$$

Conditionally on $\mathcal{F}_N^{T-1}$, the variables $(x_{i,T})$ are i.i.d. Therefore,

$$\mathbb{E}_T\bigg(\exp\bigg(\lambda_T b_N^2 Q_N^T - \frac{\lambda_T^2 b_N^2 \sigma_{N,T}^2(\psi_T)}{2}\bigg)\bigg|\mathcal{F}_N^{T-1}\bigg)$$

$$= \mathbb{E}_T\bigg(\exp\bigg(\frac{\lambda_T b_N}{\sqrt{N}}(\psi_T(x_{1,T}) - m_{N,T}(\psi_T)) - \frac{\lambda_T^2 b_N^2 \sigma_{N,T}^2(\psi_T)}{2N}\bigg)\bigg|\mathcal{F}_N^{T-1}\bigg)^N$$

$$= \mathbb{E}_T\bigg(1 + \frac{\lambda_T b_N}{\sqrt{N}}(\psi_T(x_{1,T}) - m_{N,T}(\psi_T)) - \frac{\lambda_T^2 b_N^2}{2N}\sigma_{N,T}^2(\psi_T)$$

$$+ \frac{\lambda_T^2 b_N^2}{2N}(\psi_T(x_{1,T}) - m_{N,T}(\psi_T))^2 + O\bigg(\frac{\lambda_T^3 b_N^3}{N^{3/2}}\bigg)\bigg|\mathcal{F}_N^{T-1}\bigg)^N$$

$$= (1 + \mathbb{E}_T(O(\lambda_T^3 b_N^3/N^{3/2})|\mathcal{F}_N^{T-1}))^N.$$

As $\psi_T$ is bounded, $O(\lambda_T^3 b_N^3/N^{3/2}) \leq K\lambda_T^3 b_N^3/N^{3/2}$, where $K$ does not depend on from $x_{1,T}$. Therefore,

(3.4)
$$\bigg(1 - K\frac{\lambda_T^3 b_N^3}{N^{3/2}}\bigg)^N \leq \mathbb{E}_T(e^{\lambda_T b_N^2 Q_N^T - \lambda_T^2 b_N^2 \sigma_{N,T}^2(\psi_T)/2}|\mathcal{F}_N^{T-1})$$

$$\leq \bigg(1 + K\frac{\lambda_T^3 b_N^3}{N^{3/2}}\bigg)^N.$$

Let us now deal with

$$\frac{\lambda_T^2 b_N^2}{2}\Delta_N = \bigg(\frac{\lambda_T^2 b_N^2}{2}(\sigma_{N,T}^2(\psi_T) - \sigma_T^2(\psi_T))\bigg).$$

Recall that

$$|\sigma_{N,T}^2(\psi_T) - \sigma_T^2(\psi_T)|$$

$$= \bigg|\int \psi_T^2(f_{T|T}^N - f_{T|T})\,d\mu$$

$$- \bigg(\int \psi_T(f_{T|T}^N + f_{T|T})\,d\mu\bigg)\bigg(\int \psi_T(f_{T|T}^N - f_{T|T})\,d\mu\bigg)\bigg|$$

$$\leq \bigg|\int \psi_T^2(f_{T|T}^N - f_{T|T})\,d\mu\bigg| + 2\|\psi_T\|_\infty\bigg|\int \psi_T(f_{T|T}^N - f_{T|T})\,d\mu\bigg|.$$

As $\mathbb{P}_T\{|\int \psi(f_{T|T}^N - f_{T|T})\,d\mu| > \varepsilon\} \leq \alpha(T)\exp(-N\epsilon^2/(\beta(T)\|\psi\|_\infty^2))$ for every bounded measurable $\psi$ by Lemma 2.1, we get

$$\limsup_{N\to\infty} \frac{1}{b_N^2}\log \mathbb{P}_T(|\sigma_{N,T}^2(\psi_T) - \sigma_T^2(\psi_T)| > \delta) = -\infty \qquad \forall \delta > 0.$$

In particular, $\tilde{R}_N^T$ and $\tilde{R}_N^T + (\lambda_T/2)\Delta_N$ are exponentially equivalent up to the speed $b_N^2$. We can now conclude

$$\limsup_{N\to\infty} \left| \frac{1}{b_N^2} \log \mathbb{E}_T \exp\left\{ b_N^2 \sum_{t=0}^{T-1} \lambda_t M_N^t(\psi_t) + b_N^2 \lambda_T \tilde{M}_N^T(\psi_T) \right\} \right.$$
$$\left. - \frac{\langle \boldsymbol{\lambda}_T, \mathbf{V}_T(\psi_{0:T}) \cdot \boldsymbol{\lambda}_T \rangle}{2} \right|$$

$$\stackrel{(a)}{=} \limsup_{N\to\infty} \left| \frac{1}{b_N^2} \log \mathbb{E}_T \exp\left\{ \frac{\lambda_T^2 b_N^2}{2} \Delta_N + b_N^2 \lambda_T \tilde{R}_N^T + b_N^2 \sum_{t=0}^{T-1} \lambda_t M_N^t(\psi_t) \right\} \right.$$
$$\left. - \frac{\langle \tilde{\boldsymbol{\lambda}}_{T-1}, \mathbf{V}_{T-1}(\psi_{0:T-2}, \rho_{T-1}) \cdot \tilde{\boldsymbol{\lambda}}_{T-1} \rangle}{2} \right|$$

$$\stackrel{(b)}{=} \limsup_{N\to\infty} \left| \frac{1}{b_N^2} \log \mathbb{E}_T \exp\left\{ b_N^2 \lambda_T \tilde{R}_N^T + b_N^2 \sum_{t=0}^{T-1} \lambda_t M_N^t(\psi_t) \right\} \right.$$
$$\left. - \frac{\langle \tilde{\boldsymbol{\lambda}}_{T-1}, \mathbf{V}_{T-1}(\psi_{0:T-2}, \rho_{T-1}) \cdot \tilde{\boldsymbol{\lambda}}_{T-1} \rangle}{2} \right|$$

$$\stackrel{(c)}{=} \limsup_{N\to\infty} \left| \frac{1}{b_N^2} \log \mathbb{E}_T \exp\{ b_N^2 \langle \tilde{\boldsymbol{\lambda}}_{T-1}, \mathbf{M}_N^{T-1}(\psi_{0:T-2}, \rho_{T-1}) \rangle \} \right.$$
$$\left. - \frac{\langle \tilde{\boldsymbol{\lambda}}_{T-1}, \mathbf{V}_{T-1}(\psi_{0:T-2}, \rho_{T-1}) \cdot \tilde{\boldsymbol{\lambda}}_{T-1} \rangle}{2} \right|$$

$$\stackrel{(d)}{=} 0,$$

where (a) comes from (3.4), (b) comes from the exponential equivalence, (c) follows from the definition of $\rho_{T-1}$ (see Remark 3.3) and (d) follows from the induction assumption. Therefore, (3.3) is proved and so is Theorem 3.1. □

**4. Extensions of the MDP.** In this section we extend the MDP to unbounded functions and we derive a functional MDP.

4.1. *The MDP for unbounded functions.* In this section we extend the MDP to unbounded functions. The main argument in the following proof is the use of a concentration property for i.i.d. random variables established by Ledoux [20]. For the sake of simplicity, we establish the MDP for $M_N^T(\psi_T)$ instead of $\mathbf{M}_N^T(\psi_{0:T})$. However, the same kind of results holds for $\mathbf{M}_N^T(\psi_{0:T})$.

Let $T \geq 1$ and assume the following stronger version of Assumption A-0:

ASSUMPTION A-1. *There exists a nonnegative constant $C_a$ such that for every $t \geq 1$, there exist functions $h_t^+, h_t^-$ for which*

$$C_a^{-1} h_t^-(x') \leq a_t(x, x') \leq C_a h_t^+(x') \qquad \forall (x, x') \in \mathbb{R}^d \times \mathbb{R}.$$

*Moreover,*

$$0 < C_a^{-1} \int h_t^-(x') b_t(x') \, d\mu(x') \leq C_a \int h_t^+(x') b_t(x') \, d\mu(x') < \infty \qquad \forall y \in \mathbb{R}^m.$$

REMARK 4.1. It is straightforward to check that Assumption A-1 yields Assumption A-0. Recall that $\hat{a}_t$ and $\hat{b}_t$ are defined in (1.3) and (1.4), then Assumption A-1 implies that

$$0 < C_a^{-1} \int h_t^- b_t \, d\mu \leq \hat{b}_t(x_{t-1}) \leq C_a \int h_t^+ b_t \, d\mu,$$

$$0 < C_a^{-2} \frac{h_t^-(x_t) b_t(x_t)}{\int h_t^+ b_t \, d\mu} \leq \hat{a}_t(x_{t-1}, x_t) \leq C_a^2 \frac{h_t^+(x_t) b_t(x_t)}{\int h_t^- b_t \, d\mu}.$$

Otherwise stated, the particle model coincides with a simple genetic model with strongly mixing $\hat{a}_t$-mutations and regular $\hat{b}_t$-selections.

Assumption A-1 enables us to introduce the following class of functions:

$$\mathcal{E}_T = \left\{ \psi : \mathbb{R}^d \to \mathbb{R}; \; \exists \beta > 0, \; \int_{\mathbb{R}^d} e^{\beta |\psi(x)|} h_T^+(x) b_T(x, y) \mu(dx) < \infty \right\},$$

$$\mathcal{E}_T^\alpha = \left\{ \psi : \mathbb{R}^d \to \mathbb{R}; \; \forall \beta > 0, \; \int_{\mathbb{R}^d} e^{\beta |\psi(x)|^{4\alpha/(1+2\alpha)}} h_T^+(x) b_T(x) \, \mu(dx) < \infty \right\}.$$

In the case where $0 < \alpha < \frac{1}{2}$, one can readily check that $0 < \frac{4\alpha}{1+2\alpha} < 1$ and $\mathcal{E}_T^\alpha$ becomes a set of functions with subexponential moments.

THEOREM 4.1. *Assume that Assumption A-1 holds.*

1. *In the case where $\psi_T \in \mathcal{E}_T$, then $V_T(\psi_T)$ is finite for every $T \geq 1$ and $M_N^T(\psi_T)$ satisfies the MDP with good rate function $I_T$.*
2. *Let $0 < \alpha < 1/2$ and fix $b_N = N^\alpha$. In the case where $\psi_T \in \mathcal{E}_T^\alpha$, then $V_T(\psi_T)$ is finite for every $T \geq 1$ and $M_N^T(\psi_T)$ satisfies the MDP with good rate function $I_T$.*

REMARK 4.2. In the case where $T = 0$, the problem reduces to the MDP for i.i.d. random variables and is well known (see, e.g., [3, 20]).

PROOF OF THEOREM 4.1. Since $f_{T|T}(x) \leq C_a^2 \frac{h_T^+(x) b_T(x)}{\int h_T^- b_T \, d\mu}$ by Assumption A-1, there exists $\beta > 0$ such that

$$\int e^{\beta |\psi_T|} f_{T|T}^N \, d\mu < \infty \qquad \left( \text{resp. } \forall \beta > 0, \; \int e^{\beta |\psi_T|^{4\alpha/(1+2\alpha)}} f_{T|T}^N \, d\mu < \infty \right)$$

whenever $\psi_T \in \mathcal{E}_T$ (resp. $\psi_T \in \mathcal{E}_T^\alpha$). Therefore, $m_T(\psi_T) = \int \psi_T f_{T|T} \, d\mu$ and $\int \psi_T^2 f_{T|T} \, d\mu$ are finite. In particular, $\sigma_T^2(\psi_T) < \infty$. Similarly,

$$|L_T \psi_T(x)| = \left| \int \psi_T a_T(x, \cdot) b_T(\cdot, y) \, d\mu \right| \leq C_a \int |\psi_T| h_T^+ b_T(\cdot, y) \, d\mu < \infty,$$

by Assumption A-1 and the function $L_T \psi_T(\cdot)$ is bounded. So is $L_T \psi_T(\cdot) - m_T(\psi_T) L_T \mathbb{1}(\cdot)$. Finally, $V_{T-1}(L_T \psi_T - m_T(\psi_T) L_T \mathbb{1}) < \infty$ by Theorem 3.1 and $V_T(\psi_T)$ is finite by (2.7). Define

$$\psi_T^c(x) = \psi_T(x) \mathbb{1}_{|\psi_T(x)| < c} \quad \text{and} \quad \bar{\psi}_T^c(x) = \psi_T(x) - \psi_T^c(x).$$

By Theorem 3.1, $M_N^T(\psi_T^c)$ satisfies the MDP with good rate function $I_{T,c}(x) = x^2/V_t(\psi_T^c)$. Let us now prove that

$$(4.1) \qquad \forall \delta > 0 \qquad \lim_{c \to \infty} \lim_{N \to \infty} \frac{1}{b_N^2} \log \mathbb{P}_T(|M_N^T(\bar{\psi}_T^c)| > \delta) = -\infty.$$

Condition (4.1) is sufficient to get an MDP for $M_N^T(\psi_T)$ since it asserts that $(M_N^T(\psi_T^c))_{c>0}$ is an exponential approximation of $M_N^T(\psi_T)$. Recall that $M_N^T = Q_N^T + R_N^T$ [see (2.4)]. Therefore, in order to prove (4.1), it is sufficient to prove that

$$(4.2) \qquad \forall \delta > 0 \qquad \lim_{c \to \infty} \lim_{N \to \infty} \frac{1}{b_N^2} \log \mathbb{P}_T(|R_N^T(\bar{\psi}_T^c)| > \delta) = -\infty,$$

$$(4.3) \qquad \forall \delta > 0 \qquad \lim_{c \to \infty} \lim_{N \to \infty} \frac{1}{b_N^2} \log \mathbb{P}_T(|Q_N^T(\bar{\psi}_T^c)| > \delta) = -\infty.$$

Let us first prove (4.2). As in Theorem 3.1, we first prove that $R_N^T(\bar{\psi}_T^c)$ and $\tilde{R}_N^T(\bar{\psi}_T^c)$ are exponentially equivalent. This result is not a direct consequence of Lemma 3.2 since $\bar{\psi}_T^c$ is not bounded. However, since $\bar{\psi}_T^c \in \mathcal{E}_T$ (resp. $\mathcal{E}_T^\alpha$) implies that $L_T \bar{\psi}_T^c - m_T(\bar{\psi}_T^c) L_T \mathbb{1}$ is bounded, one can prove the exponential equivalence as in the proof of Lemma 3.2. Now, since $\lim_{c \to \infty} V_{T-1}(L_T \bar{\psi}_T^c - m_T(\bar{\psi}_T^c) L_T \mathbb{1}) = 0$, one has

$$\lim_{c \to \infty} \limsup_{N \to \infty} \frac{1}{b_N^2} \log \mathbb{P}_T\{|\tilde{R}_N^T(\bar{\psi}_T^c)| > \delta\} = -\infty,$$

and (4.2) is proved by the exponential equivalence.

Let us now prove (4.3).

In the case where $\psi_T \in \mathcal{E}_T$, denote by

$$\beta(c) \triangleq C_a^2 \frac{\int |\bar{\psi}_T^c|^2 \, h_T^+ \, b_T \, d\mu}{\int h_T^- \, b_T \, d\mu}.$$

Then $0 < \sigma_{N,T}^2(\bar{\psi}_T^c) \leq \beta(c)$ which is deterministic and satisfies $\lim_{c \to \infty} \beta(c) = 0$.

Since the $(x_T^i)$'s are i.i.d. with law $f_{T|T}^N$ conditionally on $\mathcal{F}_N^{T-1}$, the large deviation upper bound for i.i.d. random variables yields

$$\mathbb{P}\left\{\left|\frac{1}{b_N\sqrt{N}}\sum_{i=1}^N(\bar{\psi}_T^c - m_T^N(\bar{\psi}_T^c))\right| > \delta\Big|\mathcal{F}_N^{T-1}\right\}$$

$$= \mathbb{P}\left\{\left|\frac{1}{N}\sum_{i=1}^N(\bar{\psi}_T^c - m_T^N(\bar{\psi}_T^c))\right| > \frac{b_N\delta}{\sqrt{N}}\Big|\mathcal{F}_N^{T-1}\right\}$$

$$\leq 2\exp\left(-N\Lambda_N^*\left(\frac{b_N\delta}{\sqrt{N}}\right)\right),$$

where the former inequality is valid for every $N \geq 1$ (see [14], Chapter 2) and $\Lambda_N^*$ is given by

$$\Lambda_N^*(x) = \sup_{\lambda \in \mathbb{R}}\left\{\lambda x - \ln \int e^{\lambda(\bar{\psi}_T^c - m_T^N(\bar{\psi}_T^c))} f_{T|T}^N d\mu\right\}$$

$$= \sup_{\lambda \in \mathbb{R}}\left\{\lambda x - \ln\left(1 + \frac{\lambda^2 \sigma_{N,T}^2(\bar{\psi}_T^c)}{2}\right.\right.$$

$$\left.\left. + \sum_{k=3}^\infty \frac{\lambda^k}{k!}\int(\bar{\psi}_T^c - m_T^N(\bar{\psi}_T^c))^k f_{T|T}^N d\mu\right)\right\}.$$

Since $-\ln(1+u) \geq -u$ for $u > -1$, one gets

$$\Lambda_N^*(x) \geq \sup_{\lambda \in \mathbb{R}}\left\{\lambda x - \frac{\lambda^2 \sigma_{N,T}^2(\bar{\psi}_T^c)}{2}\right.$$

$$\left. - \sum_{k=3}^\infty \frac{|\lambda|^k}{k!}\int(|\bar{\psi}_T^c| + |m_T^N(\bar{\psi}_T^c)|)^k f_{T|T}^N d\mu\right\}.$$

Let $x = \frac{b_N\delta}{\sqrt{N}}$ and choose $\lambda = \frac{b_N\delta}{\sqrt{N}\beta(c)}$, then

$$\Lambda_N^*\left(\frac{b_N\delta}{\sqrt{N}}\right)$$

$$\geq \frac{b_N^2\delta^2}{2N\beta(c)} - \frac{C_a^2}{\int h_T^- b_T\, d\mu}\sum_{k=3}^\infty \frac{|b_N\delta|^k}{|\sqrt{N}\beta(c)|^k k!}\int(|\bar{\psi}_T^c| + m_T^*(\bar{\psi}_T^c))^k h_T^+ b_T\, d\mu,$$

where $m_T^*(\bar{\psi}_T^c) = C_a^2 \frac{\int \bar{\psi}_T^c h_T^+ b_T\, d\mu}{\int h_T^- b_T\, d\mu}$. In particular,

$$\Lambda_N^*\left(\frac{b_N\delta}{\sqrt{N}}\right) \geq \frac{b_N^2\delta^2}{2N\beta(c)} - \frac{b_N^2}{2N}\Gamma(c,N),$$

where $\Gamma(c, N)$ is deterministic and $\lim_{N\to\infty} \Gamma(c, N) = 0$. Thus,

$$\mathbb{P}\left\{\left|\frac{1}{b_N\sqrt{N}}\sum_{i=1}^N (\bar{\psi}_T^c - m_T^N(\bar{\psi}_T^c))\right| > \delta\right\} \leq 2\exp\left(-\frac{b_N^2\delta^2}{2\beta(c)} + \frac{b_N^2}{2}\Gamma(c, N)\right)$$

and (4.3) is proved in the case where $\psi_T \in \mathcal{E}_T$ since $\lim_{c\to\infty} \beta(c) = 0$.

In the case where $\psi_T \in \mathcal{E}_T^\alpha$, denote by

$$\phi_1 = \bar{\psi}_T^c \mathbb{1}_{\{|\psi_T^c| \leq \sqrt{n}/b_N\}} \quad \text{and} \quad \phi_2 = \bar{\psi}_T^c - \phi_1.$$

Then

(4.4)
$$\mathbb{P}\left\{\left|\frac{1}{b_N\sqrt{N}}\sum_{i=1}^N (\bar{\psi}_T^c - m_T^N(\bar{\psi}_T^c))\right| > \delta\right\}$$
$$\leq \mathbb{P}\left\{\left|\frac{1}{b_N\sqrt{N}}\sum_{i=1}^N (\phi_1 - m_N^T(\phi_1))\right| > \frac{\delta}{2}\right\}$$
$$+ \mathbb{P}\left\{\left|\frac{1}{b_N\sqrt{N}}\sum_{i=1}^N (\phi_2 - m_T^N(\phi_2))\right| > \frac{\delta}{2}\right\}.$$

One can deal with the first part of the right-hand side of the inequality as done previously in order to obtain

(4.5)
$$\mathbb{P}\left\{\left|\frac{1}{b_N\sqrt{N}}\sum_{i=1}^N (\phi_1 - m_T^N(\phi_1))\right| > \frac{\delta}{2}\right\}$$
$$\leq 2\exp\left(-\frac{b_N^2(\delta/2)^2}{2\beta(c)} + \frac{b_N^2}{2}\Gamma(c, N)\right),$$

with $\lim_{N\to\infty} \Gamma(c, N) = 0$.

Let us now deal with the second part of the right-hand side of (4.4). Since $\phi_2 \in \mathcal{E}_T^\alpha$ and $b_N = N^\alpha$, one can prove (cf. [2]) that there exists $M > 0$ such that

(4.6) $$\limsup_{N\to\infty} \frac{1}{b_N^2} \log N\tilde{\mathbb{P}}(|\phi_2(\tilde{x}_i^N) - \tilde{\mathbb{E}}(\phi_2)| > u\sqrt{N}b_N) \leq -\frac{u^2}{M},$$

where $\tilde{\mathbb{P}}(dx) = \frac{h_T^+(x)b_T(x)\mu(dx)}{\int h_T^+ b_T \, d\mu}$ and $\tilde{x}_i^N$ is distributed according to $\tilde{\mathbb{P}}$. Denote by $\kappa = C_a^2 \frac{\int h_T^+ b_T \, d\mu}{\int h_T^- b_T \, d\mu}$. One gets

(4.7)
$$\mathbb{P}\{|\phi_2(x_i^N) - m_T^N(\phi_2)| > ub_N\sqrt{N}|\mathcal{F}_N^{T-1}\}$$
$$\leq \mathbb{P}\{|\phi_2(x_i^N) - \tilde{\mathbb{E}}(\phi_2)| > 2^{-1}ub_N\sqrt{N}|\mathcal{F}_N^{T-1}\}$$
$$+ \mathbb{P}\{|m_T^N(\phi_2) - \tilde{\mathbb{E}}(\phi_2)| > 2^{-1}ub_N\sqrt{N}|\mathcal{F}_N^{T-1}\}.$$

Since $m_T^N(|\phi_2|) \leq \kappa \tilde{\mathbb{E}}|\phi_2|$, there exists $N_0$ deterministic such that

$$\mathbb{P}\{|m_T^N(\phi_2) - \tilde{\mathbb{E}}(\phi_2)| > 2^{-1} u b_N \sqrt{N} | \mathcal{F}_N^{T-1}\} = 0 \quad \text{for } N \geq N_0.$$

Therefore, there exists $N_1$ deterministic and $M_2 > 0$ such that

$$N\mathbb{P}\{|\phi_2(x_i^N) - m_T^N(\phi_2)| > u b_N \sqrt{N} | \mathcal{F}_N^{T-1}\}$$

(4.8)

$$\leq \kappa \exp\left(-\frac{u^2 b_N^2}{M_2}\right) \quad \text{for } N \geq N_1.$$

With condition (4.8) in hand, a minor modification of the proof of Theorem 1 in [20] yields

(4.9) $\quad \forall \delta > 0 \quad \limsup_{N \to \infty} \frac{1}{b_N^2} \log \mathbb{P}\left\{\left|\frac{1}{b_N \sqrt{N}} \sum_{i=1}^N (\phi_2 - m_T^N(\phi_2))\right| > \frac{\delta}{2}\right\} = -\infty.$

Finally, (4.5) and (4.9) yield (4.3) in the case where $\psi \in \mathcal{E}_T^\alpha$.

It remains to identify the rate function. In fact, the exponential approximation procedure yields to the formula

$$\bar{I}(x) = \sup_{\varepsilon > 0} \liminf_{c \to \infty} \inf_{|z-x| \leq \varepsilon} I_{T,c}(z).$$

It is straightforward to check that $V_T(\psi_T^c) \to V_T(\psi_T)$ as $c \to \infty$ by the dominated convergence theorem. Since $I_{T,c}(x) = x^2/V_T(\psi_T^c)$, we easily get $\bar{I}(x) = x^2/V_T(\psi_T)$, which concludes the proof. $\square$

4.2. *A functional MDP.* We now state a functional version of the MDP. Consider the space of càdlàg functions $\mathbb{D} \triangleq D([0,1], \mathbb{R})$ and let $M_{N,T}^\psi : [0,1] \to \mathbb{D}$ be defined by

$$M_{N,T}^\psi(u) = \begin{cases} 0, & \text{if } 0 \leq u < 1/N, \\ \frac{1}{b_N \sqrt{N}} \sum_{i=1}^{[Nu]} (\psi(x_{i,T}) - m_T(\psi)), & \text{otherwise,} \end{cases}$$

where $[a]$ denotes the integer part of $a$. We shall establish the MDP for $M_{N,t}^\psi \in \mathbb{D}$. For the sake of simplicity, we only consider the particle profile involving the last generation of particles and we assume $\psi$ to be bounded. We will denote by $\mathcal{AC}_0$ the set of absolutely continuous functions $f$ from $[0,1]$ to $\mathbb{R}$ with $f(0) = 0$ [in particular, $f \in \mathcal{AC}_0 \Rightarrow f(v) \stackrel{\text{a.e.}}{=} \int_0^v \dot{f}(u)\, du$ with $\dot{f} \in L^1(du)$].

THEOREM 4.2. *Assume that Assumption* A-0 *holds and endow* $\mathbb{D}$ *with the supremum norm topology. Then the functional*

$$J_T(f) = \begin{cases} \int_0^1 \frac{\dot{f}^2(t)}{2\sigma_T^2(\psi)}\, dt + \frac{f^2(1)}{2}\left(\frac{1}{V_T(\psi)} - \frac{1}{\sigma_T^2(\psi)}\right), & \text{if } f \in \mathcal{AC}_0, \\ +\infty, & \text{otherwise,} \end{cases}$$

is a convex good rate function and $M_{N,T}^\psi$ satisfies the MDP in $\mathbb{D}$ with rate function $J_T$.

REMARK 4.3. Note that the contraction principle yields immediately Theorem 3.1. The unusual form of the rate function $J_T$ (from the point of view of Theorem 3.1) might be interpreted as the centering with the true filter density and the particles dependency with the last generation.

Proof of Theorem 4.2 is postponed to Appendix B.

## 5. Examples.

*A nonlinear observation model with additive noise.* Let $\{(X_t, Y_t); t \geq 0\}$ be a family of random variables recursively defined by

$$X_{t+1} = f_t(X_t, \varepsilon_t),$$
$$Y_t = g_t(X_t) + \eta_t,$$

where $f_t : \mathbb{R} \times \mathbb{R} \to \mathbb{R}$ (resp. $g_t : \mathbb{R} \to \mathbb{R}$) is a $\mathcal{B}(\mathbb{R} \times \mathbb{R})$-measurable [resp. $\mathcal{B}(\mathbb{R})$-measurable] function and $(\varepsilon_t)_{t \geq 0}$ [resp. $(\eta_t)_{t \geq 0}$] is a family of i.i.d. random variables.

PROPOSITION 5.1. *Assume that $\eta_t$ has a bounded and positive density with respect to the Lebesgue measure, then Assumption* A-0 *holds.*

PROOF. Let $b$ be the density of $\eta_t$ w.r.t. the Lebesgue measure. Using that $b_t(u, y) = b(y - g_t(u))$, we have

$$\gamma_t \leq \sup_x \int a_t(x, u) \sup_v b(v) \mu(du) = \sup_v b(v) < \infty.$$

Moreover, for all $x \in \mathbb{R}$, $L_t \mathbb{1}(x) = \int a_t(x, u) b(y - g_t(u)) \mu(du) > 0$ since $b(v) > 0$ for all $v \in \mathbb{R}$. The proof is complete. $\square$

*A stochastic volatility model.* Let $\{(X_t, Y_t); t \geq 0\}$ be a family of random variables recursively defined by

$$X_{t+1} = f_t(X_t, \varepsilon_t),$$
$$Y_t = \exp(X_t) \eta_t,$$

where $f_t : \mathbb{R} \to \mathbb{R}$ is a $\mathcal{B}(\mathbb{R})$-measurable function and $(\varepsilon_t)_{t \geq 0}$ [resp. $(\eta_t)$] is a family of i.i.d. random variables. We refer to [22] for more references and results on this model.

PROPOSITION 5.2. *Assume that $\eta_t \sim \mathcal{N}(0, 1)$, then Assumption* A-0 *holds.*

PROOF. Following the arguments used in the proof of Proposition 5.1, it is sufficient to prove that for all $y$, $\sup_x |b_t(x,y)| < \infty$. This can be achieved by noting that, since $b_t(x,y) = \frac{1}{\sqrt{2\pi}\exp(x)} \exp(-\frac{y^2}{2\exp(2x)})$, we have for all $y$, $\lim_{x\to\infty} b_t(x,y) = \lim_{x\to-\infty} b_t(x,y) = 0$. □

REMARK 5.1. In this example, Assumption A-0 holds and one has $\forall y \in \mathbb{R}^m$, $\inf_x b_t(x,y) = 0$.

*An example of unbounded functions.* Consider the previous framework, slightly modified:

$$X_{t+1} = f_t(X_t) + \varepsilon_t,$$
$$Y_t = g_t(X_t) + \eta_t,$$

in the particular case where $f_t$ and $g_t$ are bounded continuous and $\sigma_t = 1$. Assume, moreover, that $\varepsilon_t \sim \mathcal{N}(0,1)$ and that $\eta_t$ has a positive continuous density $u_t(z)$ with respect to the Lebesgue measure. In this case, the probability kernels are given by

$$a_t(x_{t-1}, x_t) = \frac{1}{\sqrt{2\pi}} \exp\left(-\frac{(x_t - f_{t-1}(x_{t-1}))^2}{2}\right),$$
$$b_t(x_t, y_t) = u_t(y_t - g_t(x_t)),$$

and Assumption A-1 is trivially satisfied. In particular, $b_t(x,y)$ is bounded for fixed $y$ uniformly in $x$. Moreover, one can choose $h_t^+$ as

$$h_t^+(x) = C_t e^{-(1/2)x^2 + M_t|x|},$$

where $C_t, M_t > 0$ are constants depending on $f_t$. In this case,

$$\mathcal{E}_t = \left\{\psi : \mathbb{R} \to \mathbb{R}; \exists \beta > 0, \int e^{\beta|\psi(x)|} e^{-(1/2)x^2 + M_t|x|} dx < \infty\right\},$$
$$\mathcal{E}_t^\alpha = \left\{\psi : \mathbb{R} \to \mathbb{R}; \forall \beta > 0, \int e^{\beta|\psi(x)|^{4\alpha/(1+2\alpha)}} e^{-(1/2)x^2 + M_t|x|} dx < \infty\right\}.$$

In particular,

$$\{\psi : \mathbb{R} \to \mathbb{R}; \exists K > 0, \ |\psi(x)| \leq K(1+x^2)\} \subset \mathcal{E}_t \quad \text{and}$$
$$\psi : x \mapsto \log^+(|x|) \in \mathcal{E}_t.$$

Moreover,

$$\{\psi : \mathbb{R} \to \mathbb{R}; \exists K > 0, \exists \theta \in (0,1), \ |\psi(x)| \leq K(1+x^{\theta(1+2\alpha)/(2\alpha)})\} \subset \mathcal{E}_t^\alpha,$$

and one gets interesting examples based on unbounded functions for which the MDP holds.

## APPENDIX A:

PROOF OF LEMMA 3.2. Let $K_N$ be defined by

$$K_N = \frac{1}{N}\sum_{i=1}^{N}(L_T\psi_T(x_{i,T-1}) - m_T(\psi_T)L_T\mathbb{1}(x_{i,T-1})).$$

We can express the difference as

$$|R_N^T - \tilde{R}_N^T| = \left|\frac{\sqrt{N}K_N}{b_N}\left\{\frac{1}{\sum_{i=1}^{N}L_T\mathbb{1}(x_{i,T-1})/N} - \frac{1}{\kappa_T}\right\}\right|.$$

Thus,

$$\mathbb{P}_T\{|R_N^T - \tilde{R}_N^T| > \delta\}$$

(A.1)
$$\leq \mathbb{P}_T\left\{|K_N| > \frac{L\sqrt{\delta}b_N}{\sqrt{N}}\right\} + \mathbb{P}_T\left\{\left|\frac{1}{\sum_{i=1}^{N}L_T\mathbb{1}(x_{i,T-1})/N} - \frac{1}{\kappa_T}\right| > \frac{\sqrt{\delta}}{L}\right\}.$$

Since $m_{T-1}(L_T\psi_T - m_T(\psi_T)L_T\mathbb{1}) = 0$, one can apply Lemma 2.1 to the fist part of the right-hand side,

$$\limsup_{N\to\infty}\frac{1}{b_N^2}\log\mathbb{P}_T\left\{|K_N| > \frac{L\sqrt{\delta}b_N}{\sqrt{N}}\right\} \leq -\frac{L^2\delta}{4\gamma_T^2\|\psi_T\|_\infty^2\beta(T)}\xrightarrow[L\to\infty]{}-\infty.$$

Let us now deal with the second part of the right-hand side of (A.1),

$$\left|\left(\frac{1}{N}\sum_{i=1}^{N}L_T\mathbb{1}(x_{i,T-1})\right)^{-1} - \kappa^{-1}(T)\right| > \frac{\sqrt{\delta}}{L}$$

$$\Longleftrightarrow \quad \frac{1}{\kappa_T}\left|\frac{1}{N}\sum_{i=1}^{N}L_T\mathbb{1}(x_{i,T-1}) - \kappa_T\right|$$

$$-\frac{\sqrt{\delta}}{L}\left(\frac{1}{N}\sum_{i=1}^{N}L_T\mathbb{1}(x_{i,T-1}) - \kappa_T\right) > \frac{\sqrt{\delta}}{L}\kappa_T$$

$$\Longrightarrow \quad \left(\frac{1}{\kappa_T} + \frac{\sqrt{\delta}}{L}\right)\left|\frac{1}{N}\sum_{i=1}^{N}L_T\mathbb{1}(x_{i,T-1}) - \kappa_T\right| > \frac{\sqrt{\delta}}{L}\kappa_T.$$

Denote by $\varepsilon(L) \triangleq \frac{\sqrt{\delta}}{L}\kappa_T(\frac{1}{\kappa_T} + \frac{\sqrt{\delta}}{L})^{-1}$ and apply Lemma 2.1 to conclude

$$\limsup_{N\to\infty}\frac{1}{b_N^2}\log\mathbb{P}_T\left\{\left|\frac{1}{N}\sum_{i=1}^{N}L_N\mathbb{1}(x_{i,T-1}) - \kappa_T\right| > \varepsilon(L)\right\} = -\infty$$

for every $L > 0$. Therefore, Lemma 3.2 is proved. $\square$

## APPENDIX B:

PROOF OF THEOREM 4.2. In order to prove the functional MDP, one may want to follow the lines of proof of Theorem 5.1.2 in [14] and to introduce the polygonal approximation $\check{M}^\psi_{N,T}$ of $M^\psi_{N,T}$. However, it turns out to be nontrivial to prove the exponential tightness of $\check{M}^\psi_{N,T}$ in $C[0,1]$ for the sup-norm. In particular, the level sets are useless since one can show that $\lim_{N\to\infty} \mathbb{P}(\check{M}^\psi_{N,T} \notin \{J \le a\}) = 1$. This issue is circumvented with the use of the following lemma whose proof can be found in [16], Lemma A.1. Denote by $\mathcal{U}$ the set of all subdivisions of $[0,1]$, that is,

$$U \in \mathcal{U} \iff U = \{0 = u_0 < u_1 < \cdots < u_m \le 1\}.$$

LEMMA B.1. *Let $(X_N(u); 0 \le u \le 1)_{N \ge 0}$ be a sequence of càdlàg processes defined on $(\Omega, F, \mathbb{P})$. Endow the space $\mathbb{D}$ with the uniform convergence topology and let $(\lambda(N))_{N \ge 0}$ be a sequence of positive numbers going to infinity. Assume that:*

(i) *For every $U = \{0 = u_0 < u_1 < \cdots < u_m \le 1\} \in \mathcal{U}$, $(X_N(u_1), \ldots, X_N(u_m))$ satisfies the LDP on $\mathbb{R}^m$ with speed $\lambda(N)$ and rate function $I^U$.*
(ii) *For every $\delta > 0$,*

$$\lim_{\varepsilon \to 0} \sup_{0 \le u \le 1} \limsup_{N \to \infty} \frac{1}{\lambda(N)} \log \mathbb{P}\left(\sup_{u \le v \le u+\varepsilon} |X_N(v) - X_N(u)| > \delta\right) = -\infty$$

*[convention: $\forall u > 1$, $X_N(u) \triangleq X_N(1)$].*

*Then $(X_N)$ satisfies the LDP in $\mathbb{D}$ with speed $\lambda(N)$ and rate function given by*

$$I(f) = \sup_{U \in \mathcal{U}} I^U((f(u_1), \ldots, f(u_m))), \qquad f \in \mathbb{D}.$$

*Moreover, the set $\{I < +\infty\}$ is a subset of the space $C[0,1]$ of continuous functions over $[0,1]$.*

As in Theorem 3.1, we will proceed by induction. Since the function $\psi$ is bounded and the first generation of particle is an i.i.d. sample from $a_0$, Mogulskii's theorem yields the functional MDP at time $T = 0$.

STEP 1. *The finite-dimensional MDP*. Recall that

$$M^\psi_{N,T}(u) = \frac{1}{b_N \sqrt{N}} \sum_{i=1}^{[Nu]} (\psi(x_{i,T}) - m_T(\psi))$$

$$= \frac{1}{b_N \sqrt{N}} \sum_{i=1}^{[Nu]} (\psi(x_{i,T}) - m_{N,T}(\psi)) + \frac{[Nu]}{b_N \sqrt{N}} (m_{N,T}(\psi) - m_T(\psi)),$$

and let $U = \{0 = u_0 < u_1 < \cdots < u_m \leq 1\} \in \mathcal{U}$. We shall first prove that the vector $(M_{N,T}^\psi(u_i))_{1 \leq i \leq m}$ satisfies the MDP. Since the map

$$(x_1, \ldots, x_m) \mapsto (x_1, x_2 - x_1, \ldots, x_m - x_{m-1})$$

is continuous and one-to-one, it is sufficient by the contraction principle to prove the MDP for

$$(M_{N,T}^\psi(u_i) - M_{N,T}^\psi(u_{i-1}))_{1 \leq i \leq m}.$$

Moreover, using the exponential equivalence proved in Lemma 3.2, one only has to prove the MDP for the family $(\tilde{M}_{N,T}^\psi(u_i) - \tilde{M}_{N,T}^\psi(u_{i-1}))_{1 \leq i \leq m}$, where

$$\tilde{M}_{N,T}^\psi(u) = \frac{1}{\sqrt{N}b_N} \sum_{i=1}^{[Nu]} (\psi(x_{i,T}) - m_{N,T}(\psi))$$

$$+ \frac{u}{\sqrt{N}b_N \kappa_T^2} \left( \sum_{i=1}^N L_T \psi(x_{i,T-1}) - m_T(\psi) \sum_{i=1}^N L_T \mathbb{1}(x_{i,T-1}) \right)$$

$$\triangleq Q_N^T(u) + \tilde{R}_N^T(u).$$

By the same arguments as in the proof of Theorem 3.1 (conditioning with respect to $\mathcal{F}_N^T$ to deal with $Q_N^T$, using the induction assumption to deal with $\tilde{R}_N^T$), one can show that the limit

$$\Lambda(\lambda_1, \ldots, \lambda_m) = \lim_{N \to \infty} \frac{1}{b_N^2} \log \mathbb{E}_T \exp \left( b_N^2 \sum_{i=1}^m \lambda_i (\tilde{M}_{N,T}(u_i) - \tilde{M}_{N,T}(u_{i-1})) \right)$$

is equal to

$$\Lambda_T(\lambda_1, \ldots, \lambda_m) = \frac{1}{2} \sum_{i=1}^m (u_i - u_{i-1}) \lambda_i^2 \sigma_T^2(\psi)$$

$$+ \frac{V_T(\psi)}{2\kappa_T^2} \left[ \sum_{i=1}^m \lambda_i (u_i - u_{i-1}) \right]^2.$$

Consequently, the MDP is proved for $(\tilde{M}_{N,T}^\psi(u_i) - \tilde{M}_{N,T}^\psi(u_{i-1}))_{1 \leq i \leq m}$ by Gärtner–Ellis' theorem and one can identify the rate function as

$$\check{J}_T^U(x_1, \ldots, x_m) = \sum_{i=1}^m \frac{x_i^2}{2(u_i - u_{i-1})\sigma_T^2(\psi)}$$

$$+ \frac{V_{T-1}(L_T \psi_T - m_T(\psi) L_T \mathbb{1})}{2\sigma_T^2(\psi) V_T(\psi) \kappa_T^2} \left[ \sum_{i=1}^m x_i \right]^2.$$

Thus, the MDP holds for $(M^\psi_{N,T}(u_i))_{1\leq i\leq m}$ with good rate function

$$J^U_T(x_1,\ldots,x_m) = \sum_{i=1}^m \frac{(x_i - x_{i-1})^2}{2(u_i - u_{i-1})\sigma^2_T(\psi)}$$
$$+ \frac{V_{T-1}(L_T\psi_T - m_T(\psi)L_T\mathbb{1})}{2\sigma^2_T(\psi)V_T(\psi)\kappa^2_T}(x_m)^2$$
$$= \sum_{i=1}^m \frac{(x_i - x_{i-1})^2}{2(u_i - u_{i-1})\sigma^2_T(\psi)} + \left(\frac{1}{V_T(\psi)} - \frac{1}{\sigma^2_T(\psi)}\right)\frac{x_m^2}{2}.$$

STEP 2. *The negligibility with respect to the sup-norm.* Let us prove now that $\forall \delta > 0$,

$$\lim_{\varepsilon \to 0} \sup_{0 \leq u \leq 1} \limsup_{n \to \infty} \frac{1}{b_N^2} \log \mathbb{P}_T\left(\sup_{u \leq v \leq u+\varepsilon} |M^\psi_{N,T}(v) - M^\psi_{N,T}(u)| > \delta\right) = -\infty.$$

By Lemma 3.2, it is sufficient to prove that $\forall \delta > 0$,

$$\lim_{\varepsilon \to 0} \sup_{0 \leq u \leq 1} \limsup_{N \to \infty} \frac{1}{b_N^2} \log \mathbb{P}_T\left(\sup_{u \leq v \leq u+\varepsilon} |Q^T_N(v) - Q^T_N(u)| > \delta\right) = -\infty,$$

$$\lim_{\varepsilon \to 0} \sup_{0 \leq u \leq 1} \limsup_{N \to \infty} \frac{1}{b_N^2} \log \mathbb{P}_T\left(\sup_{u \leq v \leq u+\varepsilon} |\tilde{R}^T_N(v) - \tilde{R}^T_N(u)| > \delta\right) = -\infty.$$

Note that the last limit is directly obtained by the induction assumption. Thus, we only have to establish the first one,

$$\mathbb{P}_T\left(\sup_{u \leq v \leq u+\varepsilon} |Q^T_N(v) - Q^T_N(u)| > \delta\right)$$

$$= \mathbb{P}_T\left(\sup_{u \leq v \leq u+\varepsilon} \left|\sum_{i=1}^{[Nv]} (\psi(x_{i,T}) - m_{N,T}(\psi))\right.\right.$$

(B.1)
$$\left.\left. - \sum_{i=1}^{[Nu]} (\psi(x_{i,T}) - m_{N,T}(\psi))\right| > \delta b_N \sqrt{N}\right)$$

$$= \mathbb{E}_T\left(\mathbb{P}_T\left(\sup_{u \leq v \leq u+\varepsilon} \left|\sum_{i=[Nu]+1}^{[Nv]} (\psi(x_{i,T}) - m_{N,T}(\psi))\right| > \delta b_N \sqrt{N}\bigg|\mathcal{F}^{T-1}_N\right)\right)$$

$$\leq \mathbb{E}_T\left(\mathbb{P}_T\left(\max_{1 \leq k \leq [N\varepsilon]+1} \left|\sum_{i=1}^k (\psi(x_{i,T}) - m_{N,T}(\psi))\right| > \delta b_N \sqrt{N}\bigg|\mathcal{F}^{T-1}_N\right)\right),$$

where the last inequality follows from the fact that $[Nv] - [Nu] \leq [N\varepsilon] + 1$ and by the stationarity conditionally on $\mathcal{F}^{T-1}_N$. By Ottaviani's inequality

(see [21], Chapter 6, Lemma 6.2),

$$\mathbb{P}_T\left(\max_{1\leq k\leq [N\varepsilon]+1}\left|\sum_{i=1}^{k}(\psi(x_{i,T})-m_{N,T}(\psi))\right|>\delta b_N\sqrt{N})|\mathcal{F}_N^{T-1}\right)$$

$$\text{(B.2)} \quad \leq \mathbb{P}_T\left(\left|\sum_{i=1}^{[N\varepsilon]+1}(\psi(x_{i,T})-m_{N,T}(\psi))\right|>\delta b_N\sqrt{N}/2|\mathcal{F}_N^{T-1}\right)$$

$$\times \left\{1-\max_{1\leq k\leq [N\varepsilon]}\mathbb{P}_T\left(\left|\sum_{i=k+1}^{[N\varepsilon]+1}(\psi(x_{i,T})-m_{N,T}(\psi))\right|\right.\right.$$

$$\left.\left.>\delta b_N\sqrt{N}/2|\mathcal{F}_N^{T-1}\right)\right\}^{-1}.$$

Let us first control the lower part of (B.2):

$$\mathbb{P}_T\left(\left|\sum_{i=k+1}^{[N\varepsilon]+1}(\psi(x_{i,T})-m_{N,T}(\psi))\right|>\delta b_N\sqrt{N}/2|\mathcal{F}_N^{T-1}\right)$$

$$\leq 4\,\frac{\mathbb{E}_T(|\sum_{i=k+1}^{[N\varepsilon]+1}(\psi(x_{i,T})-m_{N,T}(\psi))|^2|\mathcal{F}_N^{T-1})}{\delta^2 b_N^2 N}$$

$$= 4\,\frac{([N\varepsilon]-k+1)\sigma_{[N\varepsilon]-k+1,T}^2}{\delta^2 b_N^2 N}$$

$$\leq 16\frac{\|\psi\|_\infty^2}{\delta^2 b_N^2}\xrightarrow[N\to\infty]{}0.$$

In particular, there exists a deterministic $N_0$ such that for every $N\geq N_0$,

$$\max_{1\leq k\leq [N\varepsilon]}\mathbb{P}_T\left(\left|\sum_{i=k+1}^{[N\varepsilon]+1}(\psi(x_{i,T})-m_{N,T}(\psi))\right|>\delta b_N\sqrt{N}/2|\mathcal{F}_N^{T-1}\right)\leq 1/2.$$

Inequality (B.1) together with (B.2) yields

$$\mathbb{P}_T\left(\sup_{u\leq v\leq u+\varepsilon}|Q_N^T(v)-Q_N^T(u)|>\delta\right)$$

$$\leq 2\mathbb{P}_T\left(\frac{1}{b_N\sqrt{N}}\left|\sum_{i=1}^{[N\varepsilon]+1}(\psi(x_{i,T})-m_{N,T}(\psi))\right|>\frac{\delta}{2}\right).$$

Now,

$$\mathbb{P}_T\left(\frac{1}{b_N\sqrt{N}}\left|\sum_{i=1}^{[N\varepsilon]+1}(\psi(x_{i,T})-m_{N,T}(\psi))\right|>\frac{\delta}{2}\right)$$

$$\leq \mathbb{P}_T\left(\frac{1}{b_N\sqrt{N}}\left|\sum_{i=1}^{[N\varepsilon]+1}(\psi(x_{i,T})-m_T(\psi))\right|>\frac{\delta}{4}\right)$$

$$+\mathbb{P}_T\left(\frac{[N\varepsilon]}{\sqrt{N}b_N}|m_T(\psi)-m_{N,T}(\psi)|>\frac{\delta}{4}\right).$$

And

$$\limsup_{n\to\infty}\frac{1}{b_N^2}\log\mathbb{P}_T\left(\frac{1}{b_N\sqrt{N}}\left|\sum_{i=1}^{[N\varepsilon]+1}(\psi(x_{i,T})-m_{N,T}(\psi))\right|>\frac{\delta}{2}\right)$$

$$\leq \sup\left(\limsup_{n\to\infty}\frac{1}{b_N^2}\log\mathbb{P}_T\left(\frac{1}{b_N\sqrt{N}}\left|\sum_{i=1}^{[N\varepsilon]+1}(\psi(x_{i,T})-m_T(\psi))\right|>\frac{\delta}{4}\right);\right.$$

$$\left.\limsup_{n\to\infty}\frac{1}{b_N^2}\log\mathbb{P}_T\left(\frac{[N\varepsilon]}{\sqrt{N}b_N}|m_T(\psi)-m_{N,T}(\psi)|>\frac{\delta}{4}\right)\right).$$

By the first step of the proof (finite dimensional MDP),

$$\limsup_{n\to\infty}\frac{1}{b_N^2}\log\mathbb{P}_T\left(\frac{1}{b_N\sqrt{N}}\left|\sum_{i=1}^{[N\varepsilon]+1}(\psi(x_{i,T})-m_T(\psi))\right|>\frac{\delta}{4}\right)$$

$$\leq -\frac{\delta^2}{16\varepsilon\sigma_t^2(\psi)}+\left(\frac{1}{V_T(\psi)}-\frac{1}{\sigma_T^2(\psi)}\right)\frac{(\delta/4)^2}{2}\xrightarrow[\varepsilon\to 0]{}-\infty.$$

By Lemma 3.2, $\frac{[N\varepsilon]}{\sqrt{N}b_N}(m_T(\psi)-m_{N,T}(\psi))$ and $\tilde{R}_N^T(\varepsilon)$ are exponentially equivalent up to the speed $b_N^2$. Thus,

$$\limsup_{n\to\infty}\frac{1}{b_N^2}\log\mathbb{P}_T\left(\frac{[N\varepsilon]}{\sqrt{N}b_N}|m_T(\psi)-m_{N,T}(\psi)|>\frac{\delta}{4}\right)$$

$$=\limsup_{n\to\infty}\frac{1}{b_N^2}\log\mathbb{P}_T\left(|\tilde{R}_N^T(\varepsilon)|>\frac{\delta}{4}\right)$$

$$\leq -\frac{\delta^2}{16\varepsilon^2 V_{T-1}(L_T(\psi)-m_T(\psi)L_T(1))}\xrightarrow[\varepsilon\to 0]{}-\infty.$$

This yields the desired result:

$$\lim_{\varepsilon\to\infty}\sup_{0\leq u\leq 1}\limsup_{n\to\infty}\frac{1}{b_N^2}\log\mathbb{P}_T\left(\sup_{u\leq v\leq u+\varepsilon}|Q_N^T(v)-Q_N^T(u)|>\delta\right)=-\infty.$$

STEP 3. *The MDP.* By Steps 1 and 2, the assumptions of Lemma B.1 are satisfied. Therefore, $M_{N,t}^{\psi}$ satisfies the MDP in $\mathbb{D}$ endowed with the sup-norm topology with speed $b_N^2$ and good rate function

$$\hat{J}(f) = \sup_{U \in \mathcal{U}} J_T^U(f(u_1), \ldots, f(u_m)).$$

STEP 4. *Identification of the rate function.* The identification of the rate function is fairly standard and can be done as in the proof of Lemma 5.1.6 in [14]:

$$\hat{J}(f) = \sup_U J_T^U(f(u_1), \ldots, f(u_m))$$
$$= \int_0^1 \frac{\dot{f}^2(t)}{2\sigma_T^2(\psi)} dt + \frac{f^2(1)}{2}\left(\frac{1}{V_T(\psi)} - \frac{1}{\sigma_T^2(\psi)}\right).$$

The convexity of $J_F$ is now straightforward. This ends the proof of Theorem 4.2.

□

**Acknowledgment.** The authors would like to thank the anonymous referee for well-pointed remarks which led to a substantial improvement of the paper.

R. DOUC
CMAP
ECOLE POLYTECHNIQUE
91128 PALAISEAU CEDEX
FRANCE
E-MAIL: douc@cmapx.polytechnique.fr

A. GUILLIN
CEREMADE
UNIVERSITÉ PARIS DAUPHINE
75775 PARIS CEDEX 16
FRANCE
E-MAIL: guillin@ceremade.dauphine.fr

J. Najim
CNRS
Ecole Nationale Supérieure
    des Télécommunications
46 rue Barrault
75013 Paris
France
e-mail: najim@tsi.enst.fr